\documentclass[12pt]{article}
\usepackage[intlimits]{amsmath}
\usepackage{amssymb}
\usepackage{epsfig}

\makeatletter

\newcommand{\singlespacing}{\let\CS=\@currsize\renewcommand{\baselinestretch}{1}\tiny\CS}
\newcommand{\oneandahalfspacing}{\let\CS=\@currsize\renewcommand{\baselinestretch}{1.25}\tiny\CS}
\newcommand{\doublespacing}{\let\CS=\@currsize\renewcommand{\baselinestretch}{1.35}\tiny\CS}

\def\@citex[#1]#2{\if@filesw\immediate\write\@auxout{\string\citation{#2}}\fi
  \def\@citea{}\@cite{\@for\@citeb:=#2\do
    {\@citea\def\@citea{,\linebreak[0]\hskip0pt plus .2em}%
      \@ifundefined{b@\@citeb}%
    {{\bf ?}\@warning{Citation `\@citeb' on page \thepage\space undefined}}%
      \hbox{\csname b@\@citeb\endcsname}}}{#1}}

\newtheorem{theorem}{Theorem}
\newtheorem{remark}{Remark}

\newtheorem{preposition}{Preposition}
\newtheorem{lemma}{Lemma}
\newtheorem{corollary}{Corollary}

\newtheorem{rule-def}[theorem]{Rule}

\date{}

\begin{document}
\newcommand{\be}{\begin{equation}}
\newcommand{\ee}{\end{equation}}
\newcommand{\bea}{\begin{eqnarray}}
\newcommand{\eea}{\end{eqnarray}}
\newcommand{\nn}{\nonumber}
\newcommand{\lb}{\label}
\newcounter{saveeqn}
\newcommand{\alpheqn}{\setcounter{saveeqn}{\value{equation}}%
\stepcounter{saveeqn}\setcounter{equation}{0}%
\renewcommand{\theequation}{\mbox{\arabic{saveeqn}\alph{equation}}}}
\newcommand{\reseteqn}{\setcounter{equation}{\value{saveeqn}}%
\renewcommand{\theequation}{\arabic{equation}}}
\title{\bf A note on derivative dependent singular boundary value problems arising in physiology}
\author{R. K. Pandey$^{a,}$\thanks{Corresponding author Email: rkp@maths.iitkgp.ernet.in (R. K. Pandey)\newline $^{b}$Email: amitkverma02@yahoo.co.in (Amit Kumar Verma)\newline Acknowledgement: This work is supported by \textbf{Council of Scientific and Industrial Research} (CSIR), New Delhi, India.},~~Amit K. Verma$^{b}$\\
{\small\it{$^{a}$Department of Mathematics, Indian Institute of Technology, Kharagpur--$721302$, India}}\\\small\it{$^{b}$Mathematics Group, BITS Pilani, Pilani--$333031$, Rajasthan, India.}
}
\date{}
\maketitle

\noindent\textbf{Abstract}\vspace{.1cm}

In this note we establish existence of solutions of singular
boundary value problem $-(p(x)y^{\prime
}(x))^{\prime}=q(x)f(x,y,py')$ for $0< x\leq b$ and $y'(0)=0$,
$\alpha_{1}y(b)+\beta_{1}p(b)y^{\prime}(b)=\gamma_{1}$ with
$p(0)=0$ and $q(x)$ is integrable. Regions of multiple solutions have also been determined.
\vspace{.2cm}
\vspace{.2cm}
\begin{section}
{Introduction.}
\end{section}
We consider the singular boundary value problem of the type
\begin{eqnarray}
\label{1.1}My\equiv-\left(p(x)y^{\prime }(x)\right)^{\prime }&=&q(x)f(x,y,py'),~~0< x\leq b,\\
\label{1.2}y'(0)=0,&&\alpha_{1}y(b)+\beta_{1}p(b)y^{\prime}(b)=\gamma_{1},
\end{eqnarray}
where $\alpha_1>0$, $\beta_1\geq0$ and $\gamma_1$ is any finite
constant.

O'Regan \cite{DOR-JMAA-1996} presented existence results for the following non-resonant singular boundary value problem of limit circle type with derivative independent source function
\begin{eqnarray*}
\begin{gathered}
\frac{1}{p(t)}(p(t)y'(t))'+\mu q(t)y(t)=f(t,y(t)),~{\mathrm{a.e.,~on}}~[0,1],~\lim_{t\rightarrow0^{+}}p(t)y'(t)=0,~y(1)=0,
\end{gathered}
\end{eqnarray*}
where $p\in C[0,1]\cap C^{1}(0,1)$, $p>0$ on $(0,1)$, $q$ is measurable with $q>0$, a.e., on $(0,1]$ and $\int_{0}^{1}p(x)q(x)dx<\infty$. He did not assume $\int_{0}^{1}(1/p(t))dt<\infty$, and required the existence of the the following integral \[\int_{0}^{1}\frac{1}{p(s)}\left(\int_{0}^{s}p(x)q(x)dx\right)^{1/2}ds<\infty.\] For the singular differential equation (\ref{1.1}) this integral becomes
\begin{eqnarray}
\label{Limit-Integral}\int_{0}^{1}\frac{1}{p(s)}\left(\int_{0}^{s}q(x)dx\right)^{1/2}ds<\infty.
\end{eqnarray}
While in this paper, to establish existence of the problem (\ref{1.1})--(\ref{1.2}) we require
\begin{eqnarray}
\label{Necessary-Integral}\int_{0}^{b}\frac{1}{p(x)}\int_{0}^{x}q(t)dtdx<\infty.
\end{eqnarray}
In case of $p(x)=q(x)=x^{\alpha}$ the condition (\ref{Limit-Integral}) says that $0\leq \alpha<3$ while condition (\ref{Necessary-Integral}) is true for $\alpha\geq0$. Therefore even though we could not achieve uniqueness but we get existence of solutions for a wider class of functions $p(x)$.\\

Further we assume that $p(x)$ and $q(x)$ satisfy the following conditions\\

\noindent{\bf (A-1):}\\
(i) $p(x)>0$ in $(0,b]$,\\
(ii) $p\in C^1(0,b)$ and for some $r>b$,\\
(iii) $x\frac{p^{\prime}(x)}{p(x)}$ is analytic in $\{x:|x|<r\}$
with the Taylor series expansion
\[x\frac{p^{\prime}(x)}{p(x)}=b_{0}+b_{1}x+b_{2}x^{2}+b_{3}x^{3}+\cdots,\hspace{.2in}(b_0\geq0).\]

\noindent{\bf (A-2):}\\
(i) $q(x)>0$ in $(0,b]$,\\
(ii) $q\in L^1(0,b)$ and for some $r>b$,\\
(iii) $x^{2}\frac{q(x)}{p(x)}$ is analytic in $\{x:|x|<r\}$ with
the Taylor series expansion
\[x^{2}\frac{q(x)}{p(x)}=c_{0}+c_{1}x+c_{2}x^{2}+c_{3}x^{3}+\cdots.\]

Such problems arise frequently in applied sciences
(\cite{JBK}, \cite{PLC}, \cite{CS}) and also in physiological studies
(\cite{AM}, \cite{NAAMA}, \cite{AG}, \cite{DG} and references as there
in).

In case source function $f$ is independent of $y'$, existence-uniqueness results have been established by several researchers (\cite{MMCPNS}, \cite{RS}-\cite{RKP-AKV08},\cite{WFF-JAP-NA-2009}). Using monotone
iterative technique and Fourier-Bessel series expansion,
existence-uniqueness of solution of the differential equation
(\ref{1.1}) for $q(x)=p(x)=x^{\alpha}$ has been established by
Russell-Shampine (\cite{RS}, for $\alpha=1,2$) and Chawla-Shivkumar
(\cite{MMCPNS}, for $\alpha\geq1$) with boundary condition
$y'(0)=0$ and $y(1)=B$. While in (\cite{RKP1},\cite{RKP2}) monotone
iterative method and eigenfunction expansion are used to establish
such results for $q(x)=p(x)$ satisfying $(A-1)$ and boundary condition $y'(0)=0$ and $y(b)=B$. Recently Ford and Pennline (\cite{WFF-JAP-NA-2009}) considered the following boundary value problem
\begin{eqnarray*}
y''(x)+\frac{m}{x}y'(x)=f(x,y),&&x\in(0,1],\\
y'(0)=0,~~~~~Ay(1)+By'(1)=C,&&A>0,B\geq0,C\geq0
\end{eqnarray*}
and established the existence-uniqueness result when $\frac{\partial f}{\partial y}$ is continuous and $\frac{\partial f}{\partial y}\geq0$ on a closed region.

In the our earlier work \cite{RKP-AKV08}, we establish
existence-uniqueness result for singular boundary value
problem with derivative independent nonlinear forcing term, i.e.,
$f\equiv f(x,y)$ and boundary conditions $y'(0)=0$ and
$\alpha_1y(b)+\beta_1p(b)y'(b)=\gamma_1$. The purpose of the
present work is to extend our earlier work \cite{RKP-AKV08} for
singular boundary value problem with derivative dependent forcing
term $f(x,y,py')$.

In this work we consider the following approximation scheme for the singular boundary value problem (\ref{1.1})--(\ref{1.2})
\begin{eqnarray}
\label{1.3}Ly_{n+1}=F(x,y_{n},py'_{n}),\hspace{.2in}0< x\leq b,\hspace{.6in}\\
\label{1.4}y_{n+1}'(0)=0,\hspace{.2in}
\alpha_{1}y_{n+1}(b)+\beta_{1}p(b)y_{n+1}'(b)=\gamma_{1},
\end{eqnarray}
where
\begin{eqnarray}
\label{1.5}Ly&=&-\left(p(x)y'(x)\right)'-\lambda q(x)y(x),\\
\label{1.6}F(x,y,py')&=&q(x)f(x,y(x),p(x)y'(x))-\lambda q(x)y(x),
\end{eqnarray}
similar to the approximation scheme in \cite{MC-CDC-PH} for the
regular case $(p=q=1)$. If we start with initial approximation
$u_0$ and $v_0$ in $C[0,b] \cap C^{2}(0,b)$ such that $u_0\geq
v_0$ defined by
\begin{eqnarray*}
Mu_0\geq q(x)f(x,u_{0},pu'_{0}),\hspace{.2in}0< x\leq b,\hspace{.2in}\\
u_{0}'(0)=0,\hspace{.2in}
\alpha_{1}u_{0}(b)+\beta_{1}p(b)u_{0}'(b)\geq\gamma_{1},
\end{eqnarray*}
and
\begin{eqnarray*}
Mv_0\leq q(x)f(x,v_{0},pv'_{0}),\hspace{.2in}0< x\leq b,\hspace{.2in}\\
v_{0}'(0)=0,\hspace{.2in}
\alpha_{1}v_{0}(b)+\beta_{1}p(b)v_{0}'(b)\leq\gamma_{1}.
\end{eqnarray*}
Then the two sequences $u_n$ and $v_n$ satisfying
(\ref{1.3})--(\ref{1.4}) are monotonic and converge uniformly to
solutions $\tilde{u}$ and $\tilde{v}$ (say) of (\ref{1.1})--({\ref{1.2})
under certain conditions on $f(x,y,py')$ in the region
\[D_0=\left\{(x,y,py'):[0,b]\times[v_0,u_0]\times\mathbb{R}\right\}.\]
Functions $u_0(x)$ and $v_0(x)$ are called upper and lower solutions of (\ref{1.1})--(\ref{1.2}).
Any solution $z(x)$ in $D_0$ satisfy $\tilde{v}(x)\leq z(x)\leq \tilde{u}(x)$.

This work improves the result of Dunninger et al. \cite{DRD-JCK-JMAA} and Bobisud \cite{LEB} since we do not require the sign restriction, i.e., $y~f(x,y,0)>0$ for $|y|>M_0$ where $M_0>0$ is a constant. This result also generalizes a recent result due to Ford and Pennline (\cite{WFF-JAP-NA-2009}).

This paper is organized as follows: In Section \ref{CLBVP} we
state a result from \cite{RKP-AKV08} regarding nonnegativity of
solution of corresponding linear problem. In Section \ref{NLBVP}
we establish monotonicity of the sequences $u_n$ and $v_n$ using
the result of Section \ref{CLBVP}.

\section{Corresponding Linear Boundary Value Problem\label{CLBVP}}
The linear boundary value problem corresponding to the
approximation scheme (\ref{1.3})--(\ref{1.4}) is exactly same as
that of \cite{RKP-AKV08}, so we skip the details and state only
the following result
\begin{corollary}\label{Corollary(1)}
$(\cite{RKP-AKV08})$ If $y(x)$ satisfies
$-(p(x)y'(x))'-\lambda q(x)y=q(x)f(x)\geq0$, $0< x\leq b$,
$y'(0)=0$ and $\alpha_{1}y(b)+\beta_{1}p(b)y^{\prime}(b)=\gamma_{1}\geq0$ then
$y(x)\geq 0$ provided $\lambda< \lambda_{1}$.
\end{corollary}
\begin{remark}
The solution $y(x)$ of linear boundary value problem
$Ly(x)=q(x)f(x)$ with boundary conditions $y'(0)=0$,
$\alpha_{1}y(b)+\beta_{1}p(b)y'(b)=\gamma_{1}$ can be written as
sum of the solutions of $Ly=qf$, $y'(0)=0$,
$\alpha_{1}y(b)+\beta_{1}p(b)y'(b)=0$ and $Ly=0$,
$y'(0)=0$, $\alpha_{1}y(b)+\beta_{1}p(b)y'(b)=\gamma_{1}$. Thus
\[y(x)=\frac{\gamma_{1}u(x,\lambda)}{\alpha_{1}u(b,\lambda)+\beta_{1}p(b)u'(b,\lambda)}+\int_{0}^{b}q(t)f(t)G(x,t,\lambda)dt,\]
where $u(x,\lambda)$ is one of the two linearly independent solutions satisfying the boundary condition at $x=0$ given in Lemma $2$ of $\cite{RKP-AKV08}$ and $G(x,t,\lambda)$ is Green's function for $Ly=qf$ satisfying boundary conditions $y'(0)=0$, $\alpha_1y(b)+\beta_1y'(b)=0$.
\end{remark}
\section{Non-linear Boundary Value Problem\label{NLBVP}}
Let $f(x,y,py')$ satisfy the following
conditions:\\
{\bf(F1)} $f(x,y,py')$ is continuous on \[D_0=\{(x,y,py'): [0,b]\times [v_0,u_0]\times \mathbb{R}\};\]
{\bf(F2)} $\exists~K_1\equiv K_1(D_0)$ such that for all
$(x,y,v),(x,w,v)\in
D_0$, \[K_1(y-w)\leq f(x,y,v)-f(x,w,v)~{\mathrm{for}}~y\geq w;\]
{\bf(F3)} $\exists~0\leq L_1\equiv L_1(D_0)$ such that for all
$(x,y,v_{1}),
(x,y,v_{2})\in D_0$, \[|f(x,y,v_{1})-f(x,y,v_2)|\leq L_1|v_{1}-v_{2}|;\]
{\bf(F4)} $f(x,u_0,pu'_0)-f(x,v_0,pv'_0)-\lambda (u_0-v_0)\geq0$ for $0< x\leq b$;\\
{\bf(F5)} For all $(x,y,v)\in D_0$, $|f(x,y,v)|\leq \varphi(|v|)$
where $\varphi:[0,\infty)\rightarrow(0,\infty)$ is continuous and
satisfies
\[\int_{0}^{b}q(s)ds<\int_{0}^{\infty}\frac{ds}{\varphi(s)}.\]
Now we establish the following two Lemmas which facilitates us to
prove monotonicity of the sequences.
\begin{lemma}\label{Lemma(1)}
For $\lambda<0$ if $Ly\geq0$, $y'(0)=0$
and $\alpha_{1}y(b)+\beta_{1}p(b)y'(b)=\gamma_{1}\geq0$ then
\begin{eqnarray}
\label{3.7}(K_1-\lambda)y-L_1(sign~y')py'\geq0,\hspace{.2in}0< x\leq b,
\end{eqnarray}
provided
\begin{eqnarray}
\label{3.8} K_1-\lambda+\lambda L_1\int_{0}^{b}q(t)dt\geq0.
\end{eqnarray}
\end{lemma}
{\textbf{Proof.}} Since $y'(0)=0$ from
$Ly\geq0$ we get that $py'\geq0$ in $[0,b]$. Thus (\ref{3.7})
reduces to $(K_1-\lambda)y-L_1py'\geq0$. Integrating $Ly\geq0$
from $0$ to $x$ we get
\begin{eqnarray}
\label{3.9}\frac{p(x)y'(x)}{y(x)}\leq (-\lambda)\int_{0}^{b}q(t)dt,
\end{eqnarray}
and the result follows from (\ref{3.8}).
\begin{remark} In case
$K_1\leq 0$, condition $(\ref{3.8})$ also gives a bound on $L_1$,
\begin{eqnarray*}
-\lambda L_1\int_{0}^{b}\leq K_1-\lambda\leq -\lambda \Rightarrow L_1\leq\frac{1}{\int_{0}^{b}q(t)dt}.
\end{eqnarray*}
If $p=q=1$ and $b=1$ then the above condition reduces to $L_1\leq
1$, which is similar to Remark $4.2$ of $\cite{MC-CDC-PH}$.
\end{remark}
\begin{lemma}\label{Lemma(2)}
For $0<\lambda<\lambda_1$ and $\lambda < K_1$, if $Ly\geq0$, $y'(0)=0$ and $\alpha_{1}y(b)+\beta_{1}p(b)y'(b)=\gamma_{1}\geq0$, then
\begin{eqnarray}
\label{3.10}(K_1-\lambda)y-L_1(sign~y')py'\geq0,\hspace{.2in}0< x\leq b,
\end{eqnarray}
provided
\begin{eqnarray}
\label{3.11}1-\lambda\int_{0}^{b}\frac{1}{p(x)}\int_{0}^{x}q(t)dt dx>0,
\end{eqnarray}
and
\begin{eqnarray}
\label{3.12}
(K_1-\lambda)\left(1-\lambda\int_{0}^{b}\frac{1}{p(x)}\int_{0}^{x}q(t)dtdx\right)-\lambda
L_1\int_{0}^{b}q(t)dt\geq0.
\end{eqnarray}
\end{lemma}
{\textbf{Proof.}} From $Ly\geq0$ and
$y'(0)=0$ it is easy to see that
$py'\leq0$ in $[0,b]$. So we require to prove
\begin{eqnarray}
\label{3.13}(K_1-\lambda)y+L_1 py'\geq0,\hspace{.2in}0< x\leq b.
\end{eqnarray}
The solution of $Ly=qf\geq0$, can be written in terms of Green's
function as follows
\[y(x)=\frac{\gamma_{1}u(x,\lambda)}{\alpha_{1}u(b,\lambda)+\beta_{1}p(b)u'(b,\lambda)}+\int_{0}^{b}q(t)f(t)G(x,t)dt.\]
Now it is easy to see that to establish the result we require to
prove that
\begin{eqnarray}
\label{3.14}(K_1-\lambda)u(x,\lambda)+L_1
p(x)u'(x,\lambda)\geq0~{\mathrm{and}}~(K_1-\lambda)G(x,t)+L_1
p(x)G'(x,t)\geq0
\end{eqnarray}
for $0< x\leq b$. First we prove $(K_1-\lambda)u(x,\lambda)+L_1
p(x)u'(x,\lambda)\geq0$ where $u(x,\lambda)$ satisfies $Lu=0$,
$u'(0)=0$. Since $u(x,\lambda)\geq0$
for $0<\lambda<\lambda_{1}$ we have $pu'\leq0$. Integrating $Lu=0$
twice and using the fact that $pu'\leq0$ we get
\[-p(x)u'(x)\leq \frac{\lambda u(x)\int_{0}^{b}q(t)dt}{1-\lambda\int_{0}^{b}\frac{1}{p(s)}\int_{0}^{s}q(t)dtds}\]
and the result follows from (\ref{3.11})--(\ref{3.12}). Similarly
we can prove $(K_1-\lambda)G(x,t)+L_1 p(x)G'(x,t)\geq0$.
This completes the proof.
\begin{remark}
If $L_1=0$, i.e., $f$ is independent of $y'$, then Lemma $\ref{Lemma(1)}$ and
Lemma $\ref{Lemma(2)}$ follows from $(F2)$. In this case we can choose
$\lambda$ such that $\lambda\leq K_1$.
\end{remark}
\begin{lemma}\label{Lemma(3)}
If $u_n$ is an upper solution of
$(\ref{1.1})$--$(\ref{1.2})$ and $u_{n+1}$ is defined by
$(\ref{1.3})$--$(\ref{1.6})$ then $u_n\geq u_{n+1}$ for
$\lambda<\lambda_1$.
\end{lemma}
{\textbf{Proof.}} Let $w=u_n-u_{n+1}$. $w$ satisfies
\begin{eqnarray*}
-(pw')'-\lambda q
w=-(pu'_n)'-qf(x,u_n,py'_n)\geq0,\hspace{.2in}0< x\leq b,\\
w'(0)=0,\hspace{.2in}\alpha_{1}w(b)+\beta_{1}p(b)w'(b)\geq0,\hspace{.6in}
\end{eqnarray*}
and the result follows from Corollary \ref{Corollary(1)}.
\begin{preposition}\label{Proposition(1)}
Let $u_0$ be the upper solution of $(\ref{1.1})$--$(\ref{1.2})$,
$f(x,y,py')$ satisfy $(F1)$--$(F3)$ and $(\ref{3.8})$,
$(\ref{3.11})$, $(\ref{3.12})$ hold. Then the functions $u_n$
defined by $(\ref{1.3})$--$(\ref{1.6})$ are such that
for all $n\in \mathbb{N}$,\\
(i) $u_n$ is upper solution of $(\ref{1.1})$--$(\ref{1.2})$;\\
(ii) $u_n\geq u_{n+1}$.
\end{preposition}
{\textbf{Proof.}} We will prove this by induction. For $n=0$,
$u_0$ is an upper solution and from Lemma \ref{Lemma(3)} we have
$u_0\geq u_1$. Hence the claim is true is for $n=0$.

Let the claim be true for $n-1$, i.e., $u_{n-1}$ is an upper
solution of (\ref{1.1}) and $u_{n-1}\geq u_{n}$. Let
$w=u_{n-1}-u_{n}$. We have
\begin{eqnarray*}
-(pu'_n)'-q f(x,u_n,pu'_n)&=&
q\{f(x,u_{n-1},pu'_{n-1})-f(x,u_n,pu'_n)+\lambda (u_n-u_{n-1})\},\\
&\geq&q\{K_1(u_{n-1}-u_{n})-L_1|pu'_{n-1}-pu'_n|+\lambda (u_n-u_{n-1})\},\\
&\geq&q\{(K_1-\lambda
)w-L_1(sign~w')pw'\},
\end{eqnarray*}
and from Lemma \ref{Lemma(1)}-\ref{Lemma(2)} we get
\begin{eqnarray*}
-(pu'_n)'-q f(x,u_n,pu'_n)&\geq&0,\hspace{.2in}0< x\leq b.
\end{eqnarray*}
Thus $u_n$ is an upper solution. From Lemma \ref{Lemma(3)} we have
$u_n\geq u_{n+1}$. Hence the result follows.\\

Similarly for lower solutions we can easily deduce the following results.
\begin{lemma}\label{Lemma(4)}
If $v_n$ is a lower solution of $(\ref{1.1})$--$(\ref{1.2})$ and
$v_{n+1}$ is defined by $(\ref{1.3})$--$(\ref{1.6})$ then $v_n\leq
v_{n+1}$ for $\lambda<\lambda_1$.
\end{lemma}
\begin{preposition}\label{Proposition(2)}
Let $v_0$ be the lower solution of $(\ref{1.1})$--$(\ref{1.2})$,
$f(x,y,py')$ satisfy $(F1)$--$(F3)$ and $(\ref{3.8})$,
$(\ref{3.11})$, $(\ref{3.12})$ hold. Then the functions $v_n$
defined by $(\ref{1.3})$--$(\ref{1.6})$ are such that
for all $n\in \mathbb{N}$,\\
(i) $v_n$ is lower solution of $(\ref{1.1})$--$(\ref{1.2})$;\\
(ii) $v_n\leq v_{n+1}$.
\end{preposition}
\begin{preposition}\label{Proposition(3)}
Let $f(x,y,py')$ satisfy $(F1)$--$(F4)$ and $(\ref{3.8})$,
$(\ref{3.11})$, $(\ref{3.12})$ hold. Then for all $n\in
\mathbb{N}$ the functions $u_n$ and $v_n$ defined by
$(\ref{1.3})$--$(\ref{1.6})$ satisfy $v_n\leq u_n$.
\end{preposition}
{\textbf{Proof.}} Define, for all $i\in N$,
\[h_i(x)=f(x,u_i,pu'_i)-f(x,v_i,pv'_i)-\lambda (u_i-v_i),\hspace{.2in}0< x\leq b.\]
If $w_i=u_i-v_i$, then $w_i$ satisfies,
\[-(pw'_i)'-\lambda q w_i=qh_{i-1}.\]
We will use induction to prove this. Since $u_0\geq v_0$ we prove that $u_1\geq v_1$. Now $w_1$ is solution of
$Lw_1=qh_{0}\geq0$, $w'_1(0)=0$ and
$\alpha_{1}w_1(b)+\beta_{1}p(b)w'_1(b)=0$, from Corollary \ref{Corollary(1)} we have $w_1\geq0$. Let
$n\geq2$, $h_{n-2}\geq0$ and $u_{n-1}\geq v_{n-1}$ then we want
to prove that $h_{n-1}\geq0$ and $u_{n}\geq v_{n}$. For this consider
\begin{eqnarray*}
h_{n-1}&=&f(x,u_{n-1},pu'_{n-1})-f(x,v_{n-1},pv'_{n-1})-\lambda (u_{n-1}-v_{n-1}),\\
&\geq& K_1(u_{n-1}-v_{n-1})-L_1|pu'_{n-1}-pv'_{n-1}|-\lambda (u_{n-1}-v_{n-1}),\\
&\geq& (K_1-\lambda)w_{n-1}-L_1(sign~w'_{n-1})pw'_{n-1}.
\end{eqnarray*}
Since $w_{n-1}$ is a solution of $Lw_{n-1}=h_{n-2}\geq0$,
$w'_{n-1}(0)=0$,
$\alpha_{1}w_{n-1}(0)+\beta_{1}p(b)w'_{n-1}(b)=0$, hence from
Lemma \ref{Lemma(1)}-\ref{Lemma(2)} we have $h_{n-1}\geq0$. Thus $Lw_n=h_{n-1}\geq0$,
$w'_{n}(0)=0$ and
$\alpha_{1}w_{n}(b)+\beta_{1}p(b)w'_{n}(b)=0$ and from Corollary \ref{Corollary(1)}
have $w_n\geq0$ i.e. $u_n\geq v_n$. This completes the
proof.
\begin{lemma}\label{Lemma(5)}
If $f(x,y,py')$ satisfies $(F5)$ then there exists $R_0>0$ such that any solution of
\begin{eqnarray}
\label{3.15}-(py')'\geq qf(x,y,py'),\hspace{.2in}0< x\leq b,\hspace{.2in}\\
\label{3.16}y'(0)=0,\hspace{.2in}\alpha_{1}y(b)+\beta_{1}p(b)y'(b)\geq\gamma_{1},
\end{eqnarray}
with $y\in\left[v_0,u_0\right]$ for all $x\in[0,b]$ satisfies $\parallel py'\parallel_{\infty}< R_0$.
\end{lemma}
{\textbf{Proof.}} Since $y'(0)=0$, for each point $x\in(0,b)$ for which $p(x)y'(x)\neq0$ belongs to an interval
$(x_0,x]\subset(0,b)$ such that $p(x_0)y'(x_0)=0$ and $py'>0$ in $(x_0,x]$.
Integrating (\ref{3.15}) from $x_0$ to $x$ we get,
\[\int_{0}^{py'}\frac{ds}{\varphi(s)}\leq\int_{0}^{b}q(s)ds.\]
From (F5) we can choose $R_0>0$ such that
\begin{eqnarray*}
\int_{0}^{py'}\frac{ds}{\varphi(s)}\leq\int_{0}^{b}q(s)ds<\int_{0}^{R_0}\frac{ds}{\varphi(s)},
\end{eqnarray*}
which gives \[p(x)y'(x)< R_0.\] Similarly we consider interval
$[x,x_0)$ in which $py'<0$ and $p(x_0)y'(x_0)=0$ and we get \[-p(x)y'(x)< R_0,\] and the result
follows.

Similarly we can prove the following lemma with reverse inequalities.
\begin{lemma}\label{Lemma(6)}
If $f(x,y,py')$ satisfies $(F5)$ then there exists $R_0>0$ such that any solution of
\begin{eqnarray}
\label{3.17}-(py')'\leq qf(x,y,py'),\hspace{.2in}0< x\leq b,\hspace{.2in}\\
\label{3.18}y'(0)=0,\hspace{.2in}\alpha_{1}y(b)+\beta_{1}p(b)y'(b)\leq\gamma_{1},
\end{eqnarray}
with $y\in\left[v_0,u_0\right]$ for all $x\in[0,b]$ satisfies $\parallel py'\parallel_{\infty}< R_0$.
\end{lemma}
Now we establish our main result which is as follows:
\begin{theorem}\label{Theorem(1)}
Let $u_0$ and $v_0$ be upper and lower solution. Let $f(x,y,py')$
satisfy $(F1)$--$(F5)$ and $(\ref{3.8})$, $(\ref{3.11})$, $(\ref{3.12})$ hold then boundary value problem
$(\ref{1.1})$--$(\ref{1.2})$ has at least one solution in the
region $D_0$. If $\lambda<\lambda_1$ is chosen such that
$\lambda\leq K_1$, where $\lambda_1$ is first eigenvalue of
corresponding eigenvalue problem. The sequences $u_n$ and $v_n$
generated by $(\ref{1.3})$--$(\ref{1.6})$ with initial iterate
$u_0$ and $v_0$ converges monotonically and uniformly towards
solutions $\widetilde{u}(x)$ and $\widetilde{v}(x)$  of
$(\ref{1.1})$--$(\ref{1.2})$. Any solution $z(x)$ in $D_0$ must
satisfy $\widetilde{v}(x)\leq z(x) \leq \widetilde{u}(x)$.
\end{theorem}
{\textbf{Proof.}} From Lemma \ref{Lemma(1)}--\ref{Lemma(6)}, Preposition \ref{Proposition(1)}--\ref{Proposition(3)} we deduce that the two monotonic sequences $\{u_n\}$ and $\{v_n\}$ are bounded by $u_0$ and
$v_0$ and by Dini's Theorem their uniform convergence is assured. Let $\{u_n\}$ and $\{v_n\}$ converge uniformly to $\widetilde{u}$ and $\widetilde{v}$ respectively. Now $\{pu'_n\}$ and $\{pv'_n\}$ are uniformly bounded and from
\[\left|py'_{n}(x_1)-py'_{n}(x_2)\right|=\left|\int_{x_1}^{x_2}(py'_n)'dt\right|,\]
and uniform convergence of $\{y_n\}$, properties (A-1), (A-2) and (F-1), it is easy to prove that $py'_n$ is
equicontinuous. Now from Arzela-Ascoli's theorem there exist a uniform
convergent subsequence $\{py'_{n_{k}}\}$ of $\{py'_n\}$. Since
limit is unique so original sequence will also converge uniformly
to the same limit say $py'$. It is easy to see that if $y_n\to \widetilde{y}$ then $py'_n\to p\widetilde{y}'$. Thus we have $u_n\to\widetilde{u}$, $v_n\to\widetilde{v}$,
$pu'_n\to p{\widetilde{u}}'$ and $pv'_n\to p{\widetilde{v}}'$.

Let $G(x,t)$ be the Green's function for the linear boundary value
problem $Ly_n=0$ with homogeneous boundary conditions
$y'_n(0)=0$ and
$\alpha_{1}y_n(b)+\beta_1p(b)y'_n(b)=0$. Then solution of
(\ref{1.3})--(\ref{1.6}) can be written as
\[y_n=Ax^2+\int_{0}^{b}G(x,t)\{F(t,y_{n-1},py'_{n-1})+h(t)\}dt,\]
where $h(t)=2A(tp'(t)+p(t))+\lambda At^2q(t)$ and
$A=\frac{\gamma_{1}}{\alpha_{1}b^2+2\beta_{1}bp(b)}$. Uniform
convergence of $y_n$, $py'_n$ and continuity of $f(x,y,py')$
implies that $(1/q)F(x,y_n,py'_n)$ converges uniformly in $[0,b]$
and hence converges in mean in $L^2_q(0,b)$. Taking limit as
$n\rightarrow\infty$ and using Lemma 2.4 (\cite{TIT-EFE}, p. 27) in
$L^2_q(0,b)$ we  get
\[\widetilde{y}=Ax^2+\int_{0}^{b}G(x,t)\{F(t,\widetilde{y},p\widetilde{y}')+h(t)\}dt,\]
which is the solution of boundary value problem
(\ref{1.1})--(\ref{1.2}). Any solution $z(x)$ in $D_0$ plays the roll
of $u_0(x)$. Hence $z(x)\geq \widetilde{v}(x)$. Similarly one concludes that
$z(x)\leq \widetilde{u}(x)$.
\begin{remark}
The case when $\lambda=0$ corresponds to the case when
$f(x,y,py')\equiv f(x,py')$. In such cases the boundary value
problem $-(py')'=qf(x,py')$,
$y'(0)=0$ and
$\alpha_1y(b)+\beta_1p(b)y'(b)=\gamma_1$ can be reduced to an
initial value problem $-z'=qf(x,z)$ and $z(0)=0$ where $py'=z$.
From the assumptions on $p$, $q$ and $f(x,y,py')$ and Banach
contraction principle one can easily conclude existence-uniqueness
of solution of the nonlinear boundary value problem.
\end{remark}
\begin{remark}
Suppose in addition to the hypothesis of Theorem $\ref{Theorem(1)}$, $|f(x,y,py')|\leq N_0$ in $D_0$. Then lower solution $v_0$ and upper solution $u_0$ of $(\ref{1.1})$--$(\ref{1.2})$ may be obtained as solution of the following linear boundary value problems
\begin{gather*}
\label{LS}-(pv'_0)'+N_0q(x)=0,\hspace{.2in}0< x\leq b,\hspace{.4in}\\
\label{LSBC}v'_{0}(0)=0,\hspace{.2in}\alpha_{1}v_{0}(b)+\beta_{1}p(b)v'_{0}(b)=\gamma_{1},\\
\label{US}-(pu'_0)'-N_{0}q(x)=0,\hspace{.2in}0< x\leq b,\hspace{.4in}\\
\label{USBC}u'_{0}(0)=0,\hspace{.2in}\alpha_{1}u_{0}(b)+\beta_{1}p(b)u'_{0}(b)=\gamma_{1}.
\end{gather*}
\end{remark}

\end{document}